\documentclass[a4paper, 14pt, noxcolor
,draft]
{article}
\usepackage{maple2e}
 
\DefineParaStyle{Maple Output}
\DefineParaStyle{Maple Output}
\DefineCharStyle{2D Math}
\DefineCharStyle{2D Output}

\usepackage[a4paper]{hyperref}

\usepackage{fancybox}
\usepackage{xy}
\usepackage{epsfig}
\usepackage[T1]{fontenc} 
\usepackage[latin1]{inputenc}
\usepackage{amsfonts,amssymb,amsmath}
\usepackage{graphics} 
\usepackage{multido, xcolor}

\newtheorem{theo}{\sc Theorem}[section]

\newtheorem{rema}[theo]{\sc Remark}

\def\qed{ \ \hfil$\square$}

\makeindex

\def \F {{\mathbb F}}

\newcommand{\GL}{{\rm GL}}
\newcommand{\GSp}{{\rm GSp}}

\newcommand{\Q}{{\mathbb Q}}

\newcommand{\Z}{{\mathbb Z}}

\font\teneusm=eusm10 \font\seveneusm=eusm7 
\font\fiveeusm=eusm5 
\newfam\eusmfam 
\textfont\eusmfam=\teneusm 
\scriptfont\eusmfam=\seveneusm 
\scriptscriptfont\eusmfam=\fiveeusm

\def\mat #1,#2,#3,#4,{\left({#1\atop #3}{#2\atop #4}\right)}
\def\bra#1,{{\left\lbrace {#1}\right\rbrace}}

\def \diag {{\rm diag}}

\def \Sp{{\rm Sp}}
\def\tm{\tilde m}

\def\l1{\langle}

\newcommand{\B}{\left(\begin{array}{cc}}
\newcommand{\E}{\end{array}\right)}

\def \fns{{${}^{*)}$}}

\newcommand{\comm}[1]
{\fns\marginpar{$\boxed
{\hskip-6pt
{\small {\sf 
\begin{tabular} {l}
 #1
\end{tabular}
}
}
}
$
}
}
\def \?  {\comm{check ?}}

\usepackage{euscript}
\let\scr=\EuScript

\let\mathcal=\scr           
\def\ang#1,{{\left\langle {#1}\right\rangle}} 
\newcommand{\ds}{\displaystyle}

\def\qed{ \ \hfil$\square$}

\def \F {{\mathbb F}}

\font\teneusm=eusm10 \font\seveneusm=eusm7 
\font\fiveeusm=eusm5 
\newfam\eusmfam 
\textfont\eusmfam=\teneusm 
\scriptfont\eusmfam=\seveneusm 
\scriptscriptfont\eusmfam=\fiveeusm

\font\tengothic=eufm10
\font\sevengothic=eufm7
\font\fivegothic=eufm5
\newfam\Gothic
\textfont\Gothic=\tengothic
\scriptfont\Gothic=\sevengothic
\scriptscriptfont\Gothic=\fivegothic

\def\mat #1,#2,#3,#4,{\left({#1\atop #3}{#2\atop #4}\right)}
\def\bra#1,{{\left\lbrace {#1}\right\rbrace}}

\def \diag {{\rm diag}}

\def \Sp{{\rm Sp}}

\def\l1{\langle}

\let\scr=\EuScript

\let\mathcal=\scr           

\def\vin{{ {\tiny \mid }  
\kern-7.29pt 
\bigcup }}

\def\ang#1,{{\left\langle {#1}\right\rangle}} 

\def \cds{{\cdot{\dots}\cdot}}

\newcommand{\cross}{\times}
\normalsize

\newcounter{ncours}{\setcounter{ncours} {1}}

\def\diag {\mathop{\rm diag}\nolimits}

\def\Sp {\mathop{\rm Sp}\nolimits}

\def\qed{\quad\hbox{\hskip 1pt\vrule width 4pt height 6pt
          depth 1.5pt\hskip 1pt}}

\input amssym.def
\input amssym.tex

\font\teneusm=eusm10 \font\seveneusm=eusm7 
\font\fiveeusm=eusm5 
\newfam\eusmfam 
\textfont\eusmfam=\teneusm 
\scriptfont\eusmfam=\seveneusm 
\scriptscriptfont\eusmfam=\fiveeusm 
\def\scr#1{{\fam\eusmfam\relax#1}}

\def\finishproclaim{\par\rm
    \ifdim\lastskip<\medskipamount\removelastskip
     \penalty55\medskip\fi}

\def\proofof#1:{\par\medskip
   \noindent{\it Proof of {\rm #1}}}

\def\Ref[#1]{\par\smallskip\hang\indent%
  \llap{\hbox to\parindent{[#1]\hfil\enspace}}%
     \ignorespaces}
\def\Item#1{\par\smallskip
  \hang\indent\llap{\hbox to\parindent
     {#1\hfill\enspace}}\ignorespaces}
\def\ItemItem#1{\par\indent\hangindent2\parindent
\hbox to\parindent{#1\hfill\enspace}\ignorespaces}

\def\arrowsim{\smash{\mathop{\longrightarrow}
 \limits^{\lower1.5pt \hbox{$\scriptstyle\sim$}}}}

 \def\diag{{\rm diag}}

\def\GL{{\rm GL}}

\begin{document}

\title{On the numerator of the symplectic Hecke series of degree three}

\author{ Alexei  Panchishkin, Kirill Vankov\\
{\tt http://www-fourier.ujf-grenoble.fr/\~{}panchish}
\\
\small e-mail : panchish$@$mozart.ujf-grenoble.fr, 
FAX:  33 (0) 4 76 51 44 78}

\date{}

\maketitle

\begin{abstract}
We find a different  method to compute
the generating series in 
Shimura's conjecture for $Sp_3$, proved by Andrianov in 1967.
Formulas 
for the Satake spherical map for $\Sp_3$ are used.
\end{abstract}
\tableofcontents
\thispagestyle{empty}

\section{A formula for the  numerator of the generating series of Hecke operators of $Sp_3$}
We obtain the following  formula for the polynomial $P(v)$ 
of degree 6 in the numerator of the Hecke-Shimura generating series 
(see \cite{An87}, (3.3.72), p. 163; this series corresponds to $D_p(s)$ in  \cite{Sh}, p.825):

\begin{equation}
\begin{split}
P(x_0,\, x_1&,\, x_2,\, x_3,\, v)=\\
 =1&-\left({\displaystyle\frac{sym_{2,1,1}}{p}}+{\displaystyle\frac{(p^2+p+1)\, sym_{1,1,1}}{p^2}}+{\displaystyle\frac{sym_{1,1,0}}{p  }}\right) \, x_0^2\, v^2\\
   &+{\displaystyle\frac{p+1}{p^2}}\left(sym_{2,2,2}+sym_{2,2,1}+sym_{2,1,1}+sym_{1,1,1}\right) \, x_0^3\, v^3\\ 
   &-\left({\displaystyle\frac{sym_{3,2,2}}{p^2}}+{\displaystyle\frac{(p^2+p+1)\, sym_{2,2,2}}{p^3}}+{\displaystyle\frac{sym_{2,2,1}}{p^2}}\right) \, x_0^4\, v^4\\
   &+{\displaystyle\frac{sym_{3,3,3}}{p^3}}\, x_0^6\,v^6 \, .
\end{split}
\end{equation}

\smallskip
In this formula the notation $sym_{i_1,i_2,i_3}$ represents the symmetric polynomial of three variables 
$x_1$, $x_2$ and $x_3$ constructed in the following way:
$$
sym_{i_1,i_2,i_3} = \sum_{\sigma\in S_n / \,{\rm Stab}(x_1^{i_1}x_2^{i_2}x_3^{i_3})} \sigma (x_1^{i_1}x_2^{i_2}x_3^{i_3}),
$$
where the summation of permuted monomials $x_1^{i_1}x_2^{i_2}x_3^{i_3}$ is normalized using the stabilizer
${\rm Stab}(x_1^{i_1}x_2^{i_2}x_3^{i_3})$ so that all coefficients are equal to 1 and $i_1\ge i_2\ge i_3\ge 0$.  
The total degree of the polynomial is $i_1+i_2+i_3$.  Here $S_n = S_3$ is the symmetric group that acts naturally on polynomials in $n$ variables, where $n=3$ in our case. For example
\begin{equation}\nonumber
\begin{split}
sym_{0,0,0}&=1\\
sym_{1,0,0}&=x_1+x_2+x_3\\
sym_{1,1,0}&=x_1x_2+x_1x_3+x_2x_3\\
sym_{1,1,1}&=x_1x_2x_3\\
sym_{4,3,2}&=x_1^4x_2^3x_3^2+x_1^4x_2^2x_3^3+x_1^3x_2^4x_3^2+x_1^3x_2^2x_3^4+x_1^2x_2^4x_3^3+x_1^2x_2^3x_3^4 \, .
\end{split}
\end{equation}
Many computations presented in this article were performed using Maple 9.50 (IBM INTEL NT).  
Symmetric polynomials $sym_{i_1i_2i_3}$ (up to total degree 18) were computed using 
the coefficient of $t$ of the generating function 
\begin{equation}\nonumber
\begin{split}
&\prod_{\sigma\in S_3}(1+tx_{\sigma(1)}^{i_1}x_{\sigma(2)}^{i_2}x_{\sigma(3)}^{i_3})=(1+tx_1^{i_1}x_2^{i_2}x_3^{i_3})(1+tx_1^{i_1}x_2^{i_3}x_3^{i_2}) \\
&\qquad \cross (1+tx_1^{i_2}x_2^{i_1}x_3^{i_3})(1+tx_1^{i_2}x_2^{i_3}x_3^{i_1})(1+tx_1^{i_3}x_2^{i_1}x_3^{i_2})(1+tx_1^{i_3}x_2^{i_2}x_3^{i_1}) \, ,
\end{split}
\end{equation}
where $i_1=0,\dots ,6$, $i_2=0,\dots ,i_1$ and $i_3=0,\dots ,i_2$. Then it is normalized by dividing out 
its leading coefficient.

Let us state our result directly in terms of the Hecke operators for the symplectic group $\Sp_n$, 
defined at p. 142 of \cite{An87}.
Consider the group of positive symplectic similitudes
\begin{align}&
\mathrm{S}=\mathrm{S}^n=\GSp_n^+(\Q)=\{M\in \mathrm{M}_{2n}(\Q)\ |\ 
{}^tMJ_nM=\mu(M)J_n, \mu(M)>0\}. 
\end{align}
For the  Siegel  modular group 
$\Gamma=\Sp_n(\Z)$ consider the double cosets
\begin{align}&
(M)=\Gamma M\Gamma \subset 
\mathrm{S}, 
\end{align}
and the Hecke operators
\begin{align}&
\mathbf{T}(a)= \sum_{M\in \mathrm{SD}_n(a)}(M), 
\end{align}
where $M$ runs through the following integral matrices 
\begin{align}&
\mathrm{SD}_n(a)=\{\diag(d_1, \cdots, d_n;e_1, \cdots, e_n)\ |\ d_i| d_{i+1}, d_n| e_{n}, e_{i+1}|e_i, d_ie_i=a\}. 
\end{align}
Let us use the notation
\begin{align}&
\mathbf{T}(d_1, \cdots, d_n;e_1, \cdots, e_n)=(\diag(d_1, \cdots, d_n;e_1, \cdots, e_n)).
\end{align}
In particular we have the operators
(see p. 149 of \cite{An87}):
\begin{align}&
\mathbf{T}(p)= \mathbf{T}(\underbrace{1, \cdots, 1}_{n},\underbrace{p, \cdots, p}_{n}),\\ &
\mathbf{T}_i(p^2)= \mathbf{T}(\underbrace{1, \cdots, 1}_{n-i},\underbrace{p, \cdots, p}_{i},\underbrace{p^2, \cdots, p^2}_{n-i},\underbrace{p, \cdots, p}_{i}),
\mbox{ for } i=1, 2, \cdots, n.
\end{align}
Then their images by the spherical map $\Omega$ are given at p.159 of  \cite{An87}:
\begin{align}&
\Omega(\mathbf{T}(p))= x_0\prod_{i=0}^n(1+x_i)=\sum_{j=0}^nx_0s_j(x_1, x_2, \cdots, x_n),\\ &
\Omega(\mathbf{T}_i(p^2))= \sum_{a+b\le n, a\ge i}p^{b(a+b+1)}\mathrm{sm}_p(a-i, a)x_0^2\omega(\pi_{a, b}).
\end{align}
Here 
$\pi_{a, b}=\left(\begin{pmatrix}E_{n-a-b} & &\\ & pE_{a}&\\ & &p^2E_{b}\end{pmatrix}\right)$ is a Hecke operator for $\GL_n$, and
the coefficient $\mathrm{sm}_p(r, a)$ denotes the number of symmetric matrices of rank $r$ and order $a$ over the field $\F_p$.
This coefficient is evaluated at p.205 of  \cite{An87}:
\begin{align}\label{phi}&
\mathrm{sm}_p(r, a)=\mathrm{sm}_p(r, r)\frac{\phi_a(p)}{\phi_{r}(p)\phi_{a-r}(p)}, 
\\ & \nonumber
\mbox{ with  } \phi_{r}(x)=(x-1)(x^2-1)\cds(x^r-1) \mbox{ for   } r\ge 1, 
\mbox{ and  } \phi_{0}(x)=1.
\end{align}
In particular, we have in  the case $n=3$ that 
\begin{align}\label{Tp}&  \nonumber
\Omega(\mathbf{T}(p)) = \mathit{x}_0\,({\mathit{sym}_{1, \,1, \,0}} + 
{\mathit{sym}_{1, \,0, \,0}} + {\mathit{sym}_{1, \,1, \,1}} + 1),
\\ & \nonumber
\Omega(\mathbf{T}_1(p^2))= {\displaystyle \frac {\mathit{x}_0^{2}\,(p^{
2} - 1)\,{\mathit{sym}_{2, \,1, \,1}}}{p^{3}}}  + {\displaystyle 
\frac {\mathit{x}_0^{2}\,{\mathit{sym}_{2, \,2, \,1}}}{p}}  + 
{\displaystyle \frac {\mathit{x}_0^{2}\,{\mathit{sym}_{2, \,1, \,0
}}}{p}}  \\ &  \nonumber
\mbox{} + 
{\displaystyle \frac {\mathit{x}_0^{2}
(p-1)(3p^2+2 p+1)
{\mathit{sym}_{1, \,1, \,1}}}{p^{4}}} 
 + {\displaystyle \frac {\mathit{x}_0^{2}\,(p^{2} - 1)\,{\mathit{
sym}_{1, \,1, \,0}}}{p^{3}}}  + {\displaystyle \frac {\mathit{x}_0
^{2}\,{\mathit{sym}_{1, \,0, \,0}}}{p}},  \nonumber
\\ & 
\Omega(\mathbf{T}_2(p^2))= p^{0}\mathrm{sm}_p(0, 2)x_0^2\omega(\pi_{2, 0})
\\ &  \nonumber
+p^{4}\mathrm{sm}_p(0, 2)x_0^2\omega(\pi_{2, 1})+p^{0}\mathrm{sm}_p(1, 3)x_0^2\omega(\pi_{3, 0})
\\ & \nonumber
=x_0^2\omega (t(1,p, p))+p^{4}x_0^2\omega (t(p, p, p^2))+
\mathrm{sm}_p(1, 3)x_0^2\omega (t(p, p, p))
\\ & \nonumber
=p^{-3}x_0^2 sym_{1,1,0}+p^{-3}x_0^2 sym_{2,1,1}+
p^{-6}(p-1)(p^2+p+1)x_0^2 sym_{1,1,1},  
\\ &  \nonumber
\Omega(\mathbf{T}_3(p^2))=p^{0}\mathrm{sm}_p(0, 3)x_0^2\omega(\pi_{3, 0})= p^{-6}x_0^2x_1x_2x_3=p^{-6}x_0^2 sym_{1,1,1}, 
\end{align}
because of the equality (\ref{phi}) with $a=3, r=1$ implying $\mathrm{sm}_p(1, 3)= (p-1)(p^2+p+1)$.

Consider the polynomial $Q_3(v)$ defining the spinor zeta function
$Z(s)$ of genus three:
\begin{align}&
Q_3(v)=\\ &\nonumber (1-x_0v)(1-x_0x_1v)(1-x_0x_2v)(1-x_0x_3v)\\ &\nonumber
(1-x_0x_1x_2v)(1-x_0x_1x_3v)
(1-x_0x_2x_3v)(1-x_0x_1x_2x_3v) .
\end{align}
Following the proof at  p.159 of  \cite{An87},
  threre exist   Hecke operators  $\mathbf{t}_j\in {\mathcal L}_\Z$  with ${\mathcal L}_\Z=\Z[ \mathbf{T}(p), \mathbf{T}_1(p^2), \cdots, \mathbf{T}_{n}(p^2)]$ such that 
\begin{align}&
\sum_{j=0}^{2^n}\Omega(\mathbf{t}_j)v^j = Q_n(v)=
(1-x_0v)(1-x_0x_1v)(1-x_0x_2v)\cds (1-x_0x_1x_2\cds x_nv) .
\end{align}
Let us consider the series 
$\mathbf{R}_n(v)  = \sum_{\delta=0}^{\infty} \mathbf{T}(p^\delta )v^\delta  \in {\mathcal L}_\Z[\![v]\!]$ and 
the polynomial $\mathbf{Q}_n(v)=\sum_{j=0}^{2^n}\mathbf{t}_jv^j$ over the Hecke ring ${\mathcal L}_\Z$.

It was established by A.N.Andrianov, 
that there exist polynomials  $\mathbf{P}_n(v)$, $\mathbf{Q}_n(v)$ such that 
\begin{align}&
\mathbf{R}_n(v)  = \sum_{\delta=0}^{\infty} \mathbf{T}(p^\delta )v^\delta  =
\frac{\mathbf{P}_n(v)}{\mathbf{Q}_n(v)},
\end{align}
with the above polynomial $\mathbf{Q}_n(v)=\sum_{j=0}^{2^n}\mathbf{t}_jv^j$ of degree $2^n$, and  such that
 $\ds \mathbf{P}_{n}(v)=\sum_{j=0}^{2^n-2}\mathbf{u}_jv^j$ is a polynomial
 of degree $2^n-2$ 
with the leading term $(-1)^{n-1}p^{{n(n+1)}2^{n-2}-n^2}[\mathbf{p}]^{2^{n-1}-1} v^{2^n-2}$
(as stated in  Theorem 6 at p. 451 of \cite{An70} and at p.61 of \S 1.3, \cite{An74}).
In the following theorem we denote by $[\mathbf{p}]_n=(pE_{2n})=\mathbf{T}_n(p^2)$  the  element (3.4.48) in \cite{An70}), so that 
$\Omega([\mathbf{p}]_n)= p^{-n(n+1)/2}x_0^2x_1\cds x_n$.
\begin{theo}[see also  \cite{An67}]\label{ThP} If $n=3$, there is the following explicit polynomial presentation: 
	\begin{align}&
\mathbf{R}_3(v)  = \sum_{\delta=0}^{\infty} \mathbf{T}(p^\delta )v^\delta 
=\frac{\mathbf{P}_3(v)}{\mathbf{Q}_3(v)},
\end{align}
where
\begin{align}\label{P3}&
\mathbf{P}_3(v) 
 =1-
p^2\left( \mathbf{T}_2(p^2)+
(p^2-p+1)(p^2+p+1)
[\mathbf{p}]_3\right)v^2
+(p+1)p^4\mathbf{T}(p)[\mathbf{p}]_3v^3
\\ \nonumber 
      &     \qquad 
-p^7[\mathbf{p}]_3\left( \mathbf{T}_2(p^2)+
(p^2-p+1)(p^2+p+1)
[\mathbf{p}]_3\right) v^4 
   +p^{15}[\mathbf{p}]_3^3\,v^6 \, \in {\mathcal L}_\Z[v].
\end{align}

\end{theo}
\begin{rema} It was pointed out to the first author by S.Boecherer, that
difficulties in the problem of analytic continuation of the spinor $L$-function of genus 3
could come from
the polynomial $\mathbf{P}_3(v)$.
Indeed, this is clearly indicated  by Kurokawa's paper \cite{Ku88}, Theorem 2
in the case of the Siegel-Eisenstein series of genus 3.
A similar polynomial is discussed in the paper of Maass \cite{Maa76}.
\end{rema}
\section{Explicit form of Shimura's conjecture for $Sp_3$}
Our result gives a complement to the solution  of Shimura's conjecture by A.N.Andrianov in \cite{An67} (see  also \cite{An68} and \cite{An69}).
This conjecture was stated in  \cite{Sh}, at p.825 as follows:
\begin{quote}
``In general, it is plausible that $D_p(X)=E(X)/F(X)$ with polynomials $E(X)$  
and $F(X)$ in $X$ with integral coefficients of degree $2^n-2$ and $2^n$, respectively'' 
\end{quote}
({\it i.e. with coefficients in } ${\mathcal L}_\Z=\Z[ \mathbf{T}(p), \mathbf{T}_1(p^2), \cdots, \mathbf{T}_{n}(p^2)]$).
\begin{theo}[see also  \cite{An67}] If $n=3$, one has the following explicit polynomial presentation: 
\begin{align}&
\mathbf{Q}_3(v) 
 = 1 - \mathbf{T}(p)\,v 
\\ & \nonumber
+p\big(\,
\mathbf{T}_1(p^2) 
+ (p^{2} + 1)\,\mathbf{T}_2(p^2) 
+ 
(1+p^2)^2
- p^{3}\mathbf{T}(p)\,\big( \mathbf{T}_2(p^2)
+ [\mathbf{p}]_3\big)
v^{3} 
\\ & \nonumber
+ p^6\big(
 - 2\,p\,\mathbf{T}_1(p^2)\,[\mathbf{p}]_3 + \,\mathbf{T}_2(p^2)
\mbox{}  - 2( p -1)\,\mathbf{T}_2(p^2)\,[\mathbf{p}]_3
 \\ &\nonumber
-(p^2+2p-1)(p^2-p+1)(p^2+p+1))
[\mathbf{p}]_3^{2} +
[\mathbf{p}]_3\,
\mathbf{T}(p)^{2}\big)v^{4} \\ &\nonumber
\mbox{} 
- \,p^{9}\,[\mathbf{p}]_3^{}\,\mathbf{T}(p)\,\big ( 
\mathbf{T}_2(p^2) +[\mathbf{p}]_3\big )\,
v^{5}   \\ &\nonumber
+p^{13}\,[\mathbf{p}]_3^{2}\,\big(\,\mathbf{T}_1(p^2) + (
p^{2} + 1)\,\mathbf{T}_2(p^2) + 
(p^2+1)^2
[\mathbf{p}]_3\big)\,v^{6} \\ &\nonumber
\mbox{} - \,p^{18}\,[\mathbf{p}]_3^{3}\,\mathbf{T}(p)
\,v^{7} +\,p^{24}\,[\mathbf{p}]_3^{4}\,v^{8}
\in {\mathcal L}_\Z[v].
\end{align}

\end{theo}
Proof follows the same lines as that of Theorem \ref{ThP}.
We compute an expression for $\Omega(\mathbf{Q}_3(v))$ in terms of
$sym_{i_1,i_2,i_3}$:
\begin{align}& \nonumber
\Omega(\mathbf{Q}_3(v))
=1 - \mathit{x}_0^{}\,({\mathit{sym}_{1, \,1, \,1}} + {\mathit{sym}_{1, \,1, \,0}} + {
\mathit{sym}_{1, \,0, \,0}} +  1)\,
v \\ &  \nonumber
\mbox{} + \mathit{x}_0^{2}\,(4\,{\mathit{sym}_{1, \,1, \,1}} + {
\mathit{sym}_{1, \,0, \,0}} + 2\,{\mathit{sym}_{2, \,1, \,1}} + 2
\,{\mathit{sym}_{1, \,1, \,0}} + {\mathit{sym}_{2, \,1, \,0}} + {
\mathit{sym}_{2, \,2, \,1}})\,v^{2} - \\ &  \nonumber
-\mathit{x}_0^{3}({\mathit{sym}_{3, \,1, \,1}} + {\mathit{sym}_{1, 
\,1, \,0}} + 4\,{\mathit{sym}_{2, \,2, \,1}} + 4\,{\mathit{sym}_{
1, \,1, \,1}} + {\mathit{sym}_{2, \,1, \,0}} + {\mathit{sym}_{2, 
\,2, \,0}} + 4\,{\mathit{sym}_{2, \,2, \,2}} \\ &  \nonumber
\mbox{} + {\mathit{sym}_{3, \,2, \,2}} + {\mathit{sym}_{3, \,2, 
\,1}} + 4\,{\mathit{sym}_{2, \,1, \,1}})v^{3}\mbox{} + 
\mathit{x}_0^{4}({\mathit{sym}_{3, \,1, \,1}} + {\mathit{sym}_{1, \,1, \,1}}
 + {\mathit{sym}_{3, \,3, \,1}} \\ &  
\mbox{} + {\mathit{sym}_{4, \,2, \,2}} + 2\,{\mathit{sym}_{2, \,1
, \,1}} + 4\,{\mathit{sym}_{3, \,2, \,2}} + 2\,{\mathit{sym}_{3, 
\,2, \,1}} + {\mathit{sym}_{2, \,2, \,0}} + 8\,{\mathit{sym}_{2, 
\,2, \,2}} \\ &  \nonumber
\mbox{} + 2\,{\mathit{sym}_{3, \,3, \,2}} + {\mathit{sym}_{3, \,3
, \,3}} + 4\,{\mathit{sym}_{2, \,2, \,1}})v^{4}\mbox{} - 
\mathit{x}_0^{5}({\mathit{sym}_{4, \,3, \,3}} + {\mathit{sym}_{4, \,3, \,2
}} + {\mathit{sym}_{2, \,2, \,1}} \\ &  \nonumber
\mbox{} + 4\,{\mathit{sym}_{2, \,2, \,2}} + 4\,{\mathit{sym}_{3, 
\,3, \,2}} + {\mathit{sym}_{3, \,3, \,1}} + 4\,{\mathit{sym}_{3, 
\,2, \,2}} + 4\,{\mathit{sym}_{3, \,3, \,3}} + {\mathit{sym}_{4, 
\,2, \,2}} + {\mathit{sym}_{3, \,2, \,1}} 
)v^{5}
 \\ &  \nonumber
 + \mathit{x}_0^{6}\,(2\,{\mathit{sym}_{3, \,3, \,2}}
 + {\mathit{sym}_{3, \,2, \,2}} + 2\,{\mathit{sym}_{4, \,3, \,3}}
 + 4\,{\mathit{sym}_{3, \,3, \,3}} + {\mathit{sym}_{4, \,3, \,2}}
 + {\mathit{sym}_{4, \,4, \,3}})\,v^{6} \\ &  \nonumber
\mbox{} - \mathit{x}_0^{7}\,({\mathit{sym}_{4, \,3, \,3}} + {
\mathit{sym}_{3, \,3, \,3}} + {\mathit{sym}_{4, \,4, \,4}} + {
\mathit{sym}_{4, \,4, \,3}})\,v^{7} + \mathit{x}_0^{8}\,{\mathit{
sym}_{4, \,4, \,4}}\,v^{8}.
 \end{align}
Then we use the polynomial expressions (\ref{Tp}) 
for the generators of the Hecke ring.
Using these generators,  we may look for a solution 
in the following form:
\begin{align}&
\Omega(\mathbf{Q}_3(v))
= 1 - \Omega(\mathbf{T}(p))v
 + 
\big(\mathit{K}_{T1p2}\Omega(\mathbf{T}_1(p^2)) + \mathit{K}_{T2p2}\Omega(\mathbf{T}_2(p^2))
 \\ & \nonumber
 + \mathit{K}_{T3p2}\Omega([\mathbf{p}]_3) + \mathit{K}_{TpTp}\Omega(\mathbf{T}(p))^{2}\big)\,v^{2} 
  \\ & \nonumber
+\big(\mathit{K}_{TpT1p2}\Omega(\mathbf{T}(p)\mathbf{T}_1(p^2)) 
+ 
\mathit{K}_{TpT2p2}\Omega(\mathbf{T}(p)\mathbf{T}_2(p^2)) 
 \\ & \nonumber
+ \mathit{K}_{TpT3p2}\,
\Omega(\mathbf{T}(p)[\mathbf{p}]_3) + \mathit{K}_{TpTpTp}\Omega(\mathbf{T}(p))^{3}\big)
v^{3}\mbox{}
 \\ & \nonumber
 + 
\big(\mathit{K}_{T1p2T1p2}\Omega(\mathbf{T}_1(p^2))^{2} 
+ \mathit{K}_{T1p2T2p2}\Omega(\mathbf{T}_1(p^2)\mathbf{T}_2(p^2)) 
\\ & \nonumber
+ 
\mathit{K}_{T1p2T3p2}\,\Omega(\mathbf{T}_1(p^2)[\mathbf{p}]_3)
\mbox{} + \mathit{K}_{T2p2T2p2}\Omega(\mathbf{T}_2(p^2))^{2}
\\ & \nonumber
 + 
\mathit{K}_{T2p2T3p2}\,\Omega(\mathbf{T}_2(p^2)[\mathbf{p}]_3) + \mathit{K}_{T3p2T3p2}\,
\Omega([\mathbf{p}]_3^{2}) \\ &\nonumber
\mbox{} + \mathit{K}_{T1p2TpTp}\,\Omega(\mathbf{T}_1(p^2)\mathbf{T}(p)^{2}) + 
\mathit{K}_{T2p2TpTp}\Omega(\mathbf{T}_2(p^2)\mathbf{T}(p)^{2}) 
\\ &\nonumber
+ \mathit{K}_{T3p2TpTp}\,\Omega([\mathbf{p}]_3\mathbf{T}(p)^{2}) 
\mbox{} + \mathit{K}_{TpTpTpTp}\,\Omega(\mathbf{T}(p))^{4}\big)v^{4}\\
&\nonumber
+
\Omega(\mathbf{t}_5)v^5+\Omega(\mathbf{t}_6)v^6
+\Omega(\mathbf{t}_7)v^7+\,p^{24}\,\Omega([\mathbf{p}]_3)^{4}v^8.  
\end{align} 

It is not too difficult to resolve the resulting equations in the indeterminate coefficients:\\
$
K_{TpT1p2} = 0, K_{TpTpTp} = 0, K_{TpT2p2} = -p^3, K_{TpT3p2} =
-p^3,
K_{T2p2} = p^3+p$, \\ $K_{T3p2} = 
p(1+p^2)^2,
K_{T1p2} = p, K_{TpTp} = 0$,
$
K_{T2p2TpTp}=0$, \\ 
$K_{T1p2T3p2}=-2p^7,
K_{T2p2T3p2}= -2p^7+2p^6, 
K_{T1p2T1p2}= 0,
K_{T1p2T2p2}= 0$, \\ 
$K_{T2p2T2p2}= p^6$, 
$K_{T1p2TpTp}= 0, 
K_{T3p2TpTp}=p^6, 
K_{TpTpTpTp}= 0$,\\ $
K_{T3p2T3p2}=
-p^6(p^2+2p-1)(p^2-p+1)(p^2+p+1)$.

Then we find the remaining coefficients $\mathbf{t}_5$, $\mathbf{t}_6$, $\mathbf{t}_7$ using the functional equation \cite{An87}, p.164
(3.3.79):   
$\mathbf{t}_{8-i}=(p^6[\mathbf{p}]_3)^{4-i}\mathbf{t}_i$ ($i=0, \cdots, 8$), 
compare with formulas in 
\cite{An67} and in 
\cite{Evd}. \qed

\section{An identity involving $\omega (t(1,p^{\lambda_2}, p^{\lambda_3}))$ and $Q_3(v)$}
The theory of Hecke rings for the symplectic group is developed in \cite{Sh}, 
\cite{An87} and \cite{AnZh95} (Ch. 3).  
At page 150 of \cite{AnZh95} we have the following identity for the spherical map 
\begin{equation}\nonumber
\begin{split}
R_n(v) & = \sum_{\delta=0}^{\infty} \Omega (\mathbf{T}^n(p^\delta ))v^\delta \\
 & = \sum_{\delta=0}^{\infty}\quad\sum_{0\le\delta_1\le\cdots\le\delta_n\le\delta}p^{n\delta_1+(n-1)\delta_2+\cdots+\delta_n}\omega(t(p^{\delta_1},\cdots,p^{\delta_n}))(x_0v)^\delta \, .
\end{split}
\end{equation}

This identity for formal generating series of elements of Hecke ring allows 
to reduce computations in the local Hecke rins of the symplectic group 
to computations in polynomial rings by applying the spherical map $\Omega$ 
to elements $\mathbf{T}^n(p^\delta ) \in L_p^n(q)$ of Hecke ring for the symplectic 
group or spherical map $\omega$ to elements  
$t(p^{\delta_1},\cdots,p^{\delta_n}) = \sum_j a_j(\Lambda g_j) \in L_{\Q}(\Lambda ^n,G_p^n)$ 
of Hecke ring for the general linear group. 
Detailed definition of spherical maps 
for Hecke elements as well as definitions and a structure of left cosets $\Lambda g_j$ 
and Hecke pairs $\Lambda ^n,G_p^n$ can be found in \cite{AnZh95} chapter 3 paragraphs 
2 and 3, see also \cite{An87} chapter 3.  
For generating elements $\pi_i(p) = \pi_i^n(p) = (\diag(1,\dots ,1,p,\dots ,p))$ with 1 
on the diagonal listed $(n-i)$ times and then $p$ listed $i$ times ($1\le i \le n$) 
the images elements under the map $\omega = \omega_p^n$ are given by the formulas
\begin{equation}\nonumber
\omega (\pi_i^n(p))=p^{- \langle i \rangle }s_i(x_1,\dots ,x_n) \quad (1\le i \le n),
\end{equation}
where
\begin{equation}\nonumber
s_i(x_1,\dots ,x_n) = \sum_{1\le \alpha_1 < \cdots < \alpha_i \le n} x_{\alpha_1} \cdots x_{\alpha_i}
\end{equation}
is the $i$th elementary symmetric polynomial (different then previously defined $sym_{i_1,i_2,i_3}$).
For an arbitrary element $t$ the map is defined by
\begin{equation}\nonumber
\omega(t)=\sum_j a_j\omega ((\Lambda g_j)) \, .
\end{equation}

Examples of computation for cases $n=1$ and $n=2$ are given in \cite{Sh}, at page 824, and in the book \cite{AnZh95}, at page 150:
\begin{align}&
R_1(v)=[(1-x_0v)(1-x_0x_1v)]^{-1}\, ,\\
&
R_2(v)=[(1-x_0v)(1-x_0x_1v)(1-x_0x_2v)(1-x_0x_1x_2v)]^{-1}\left( 1-{\displaystyle \frac{x_0^2x_1x_2v^2}{p}}\right) \, .
\end{align}
In this paper we consider $n=3$.  Acting analogously to the case $n=2$ let 
\begin{equation}
\begin{split}
\delta_2 &=\delta_1+\delta_1^\prime\\
\delta_3 &=\delta_1+\delta_2^\prime\\
\delta &=\delta_1+\delta^\prime\\
\delta^\prime &=\delta_2^\prime+\beta
\end{split}
\end{equation}
where $0\le\delta_1^\prime\le\delta_2^\prime\le\delta^\prime$, $\beta\ge0$.  Then
\begin{equation}\nonumber
\begin{split}
&R_3(v)=\sum_{\delta =0}^\infty \, \sum_{0\le \delta_1 \le \delta_2 \le \delta_3 \le \delta} p^{3\delta_1 +2\delta_2 +\delta_3} \omega (t(p^{\delta_1},p^{\delta_2},p^{\delta_3})) (x_0 v)^\delta\\
&=\sum_{\delta_1\ge 0}\sum_{\beta\ge 0, 0\le\delta_1^\prime\le\delta_2^\prime} (x_0v)^{\delta_1+\delta_2^\prime+\beta} \left(\frac{x_1x_2x_3}{p^6}\right)^{\delta_1} \omega(t(1,p^{\delta_1^\prime},p^{\delta_2^\prime})) p^{3\delta_1+2(\delta_1+\delta_1^\prime)+(\delta_1+\delta_2^\prime)}\\
&=\sum_{\delta_1\ge 0}\sum_{\beta\ge 0, 0\le\delta_1^\prime\le\delta_2^\prime} (x_0v)^{\delta_1+\delta_2^\prime+\beta} \left(\frac{x_1x_2x_3}{p^6}\right)^{\delta_1} p^{6\delta_1+2\delta_1^\prime+\delta_2^\prime} \omega(t(1,p^{\delta_1^\prime},p^{\delta_2^\prime}))\\
&=\sum_{\delta_1\ge 0}\sum_{\beta\ge 0} (x_0vx_1x_2x_3)^{\delta_1} (x_0v)^\beta \sum_{0\le\delta_1^\prime\le\delta_2^\prime} 
\omega(t(1,p^{\delta_1^\prime},p^{\delta_2^\prime})) p^{2\delta_1^\prime+\delta_2^\prime} (x_0v)^{\delta_2^\prime}\\
&=[(1-x_0v)(1-x_0x_1x_2x_3v)]^{-1} \sum_{0\le\delta_1^\prime\le\delta_2^\prime} \omega(t(1,p^{\delta_1^\prime},p^{\delta_2^\prime})) 
p^{2\delta_1^\prime+\delta_2^\prime} (x_0v)^{\delta_2^\prime}\\
&=[(1-x_0v)(1-x_0x_1v)(1-x_0x_2v)(1-x_0x_3v)(1-x_0x_1x_2v)(1-x_0x_1x_3v)\\
&\ \ \ \ \cross(1-x_0x_2x_3v)(1-x_0x_1x_2x_3v)]^{-1} P_3(v)
\end{split}
\end{equation}
Here $P_3(v)$ denotes a polynomial of degree 6 as stated in the Theorem 6 (page 451) of \cite{An70}.  
This rational polynomial presentation is proved in \cite{An69} for Hecke series and $\zeta$-functions
of the groups $GL_n$ and $SP_n$ over local fields.  Further theory and applications were developed 
for genus 2 in the work \cite{An74}.

It follows that
\begin{align}&
P_3(v)= \sum_{0\le\delta_1^\prime\le\delta_2^\prime} \omega(t(1,p^{\delta_1^\prime},p^{\delta_2^\prime})) p^{2\delta_1^\prime+\delta_2^\prime} (x_0v)^{\delta_2^\prime}\times \\
&\nonumber \times [(1-x_0x_1v)(1-x_0x_2v)(1-x_0x_3v)(1-x_0x_1x_2v)(1-x_0x_1x_3v)(1-x_0x_2x_3v)],
\end{align}
and our main result follows from the explicit computation of  the coefficients $\omega(t(1,p^{\delta_1^\prime},p^{\delta_2^\prime}))$
given in the next section.

\section{A computation for the images of the Hecke operators under the spherical map}
The formula for $P_3$ is obtained using the following computation for the images
$\omega ({t}(1,p^{\lambda_2},p^{\lambda_3}))$ of the Hecke operators under the spherical 
map for the group $\Lambda=GL_3(\Z)$.  
Note that the notation $\Omega$ used in the article \cite{An70}
corresponds to  $\omega$ in our formulas by the substitution of 
$x_1$ by $x_1/p$, $x_2$ by $x_2/p$ and $x_3$ by $x_3/p$.  
We used formula (1.7) for $\Omega$ at 
page 432 of \cite{An70} and adopted it for $\omega$.  
In the notations of that article we have that
$W$ is the group $S_n = S_3$, the set $\Sigma = \{(1,-1,0),(1,0,-1),(0,1,-1)\}$, $q=p$ and 
$\lambda = (0,\delta_1^{'},\delta_2^{'})$.  
The expression for the polynomial $c(x)$ from \cite{An70}
takes the following form
$$c(x_1,x_2,x_3)=\frac{(x_2-x_1/p)(x_3-x_1/p)(x_3-x_2/p)}{(x_2-x_1)(x_3-x_1)(x_3-x_2)}\,.$$
In the case of $n=3$ we have 28 different possibilities for $\omega(t(1,p^{\lambda_2},p^{\lambda_3}))$,
which we only need in order to compute the polynomial $P_3(v)$ of degree 6:

\begin{itemize}
\item[1)]
 $\omega(t(1,1,1))
 = 1$
\item[2)]
 $\omega(t(1,1,p))
 = {\displaystyle \frac{\mathit{sym}_{1,0,0}}{p}}$
\item[3)]
 $\omega(t(1,p,p))
 = {\displaystyle \frac{\mathit{sym}_{1,1,0}}{p^3}}$
\item[4)]
 $\omega(t(1,1,p^2))
 = {\displaystyle \frac{(p-1)\,\mathit{sym}_{1,1,0}}{p^3}}
 + {\displaystyle \frac{\mathit{sym}_{2,0,0}}{p^2}}$
\item[5)]
 $\omega(t(1,p,p^2))
 = {\displaystyle \frac{(2p^2-p-1)\,\mathit{sym}_{1,1,1}}{p^6}}
 + {\displaystyle \frac{\mathit{sym}_{2,1,0}}{p^4}}$
\item[6)]
 $\omega(t(1,p^2,p^2))
 = {\displaystyle \frac{(p-1)\,{\mathit{sym}_{2,1,1}}}{p^7}}
 + {\displaystyle \frac{\mathit{sym}_{2,2,0}}{p^6}}$
\item[7)]
 $\omega ({t}(1,1,p^{3}))
 = {\displaystyle \frac {(p^{2}  - 2\,p + 1)\,{\mathit{sym}_{1,1,1}}}{p^{5}}}
 + {\displaystyle \frac {(p^{2}-p)\,{\mathit{sym}_{2,1,0}}}{p^{5}}}
 + {\displaystyle \frac {{\mathit{sym}_{3,0,0}}}{p^{3}}}$
\item[8)]
 $\omega ({t}(1,p,p^{3}))
 = {\displaystyle \frac {(2\,p-2)\,{\mathit{sym}_{2,1,1}}}{p^{6}}}
 + {\displaystyle \frac {(p - 1)\,{\mathit{sym}_{2,2,0}}}{p^{6}}}
 + {\displaystyle \frac {{\mathit{sym}_{3,1,0}}}{p^{5}}}$
\item[9)]
 $\omega ({t}(1,p^{2},p^{3}))
 = {\displaystyle \frac {(2\,p-2)\,{\mathit{sym}_{2,2,1}}}{p^{8}}}
 + {\displaystyle \frac {(p - 1)\,{\mathit{sym}_{3,1,1}}}{p^{8}}} 
 + {\displaystyle \frac {{\mathit{sym}_{3,2,0}}}{p^{7}}}$ 
\item[10)]
 $\omega ({t}(1,p^{3},p^{3}))
 = {\displaystyle \frac {(p^{2} -2\,p + 1)\,{\mathit{sym}_{2,2,2}}}{p^{11}}} 
 + {\displaystyle \frac {(p^{2}-p)\,{\mathit{sym}_{3,2,1}}}{p^{11}}} 
 + {\displaystyle \frac {{\mathit{sym}_{3,3,0}}}{p^{9}}}$
\item[11)]
 $\omega ({t}(1,1,p^{4}))
 = {\displaystyle \frac {(p^2 - 2\,p + 1)\,{\mathit{sym}_{2,1,1}}}{p^{6}}} 
 + {\displaystyle \frac {(p^{2}-p)\,{\mathit{sym}_{2,2,0}}}{p^{6}}} 
 + {\displaystyle \frac {(p^{2}-p)\,{\mathit{sym}_{3,1,0}}}{p^{6}}} \\
 + {\displaystyle \frac {{\mathit{sym}_{4,0,0}}}{p^{4}}}$
\item[12)]
 $\omega ({t}(1,p,p^{4}))
 = {\displaystyle \frac {(2\,p^{2} - 3\,p+1)\,{\mathit{sym}_{2,2,1}}}{p^{8}}} 
 + {\displaystyle \frac {(2\,p^{2} - 2\,p)\,{\mathit{sym}_{3,1,1}}}{p^{8}}} 
 + {\displaystyle  \frac {(p^{2}-p)\,{\mathit{sym}_{3,2,0}}}{p^{8}}} \\
 + {\displaystyle  \frac {{\mathit{sym}_{4, \,1, \,0}}}{p^{6}}} $
\item[13)]
 $\omega ({t}(1,p^{2},p^{4}))
 = {\displaystyle \frac {( - 4\,p^{2} + 3\,p^{3} + 2\,p-1)\,{\mathit{sym}_{2,2,2}}}{p^{11}}} 
 + {\displaystyle \frac {(2\,p^3 - 3\,p^{2} + p)\,{\mathit{sym}_{3,2,1}}}{p^{11}}} \\	
 + {\displaystyle \frac {(p^{3} - p^{2})\,{\mathit{sym}_{3,3,0}}}{p^{11}}} 
 + {\displaystyle \frac {(p^{3} - p^{2})\,{\mathit{sym}_{4,1,1}}}{p^{11}}} 
 + {\displaystyle \frac {{\mathit{sym}_{4,2,0}}}{p^{8}}}$
\item[14)]
 $\omega ({t}(1,p^{3},p^{4}))
 = {\displaystyle \frac {(2\,p^{2} - 3\,p+ 1)\,{\mathit{sym}_{3,2,2}}}{p^{12}}} 
 + {\displaystyle \frac {(2\,p^{2} - 2\,p)\,{\mathit{sym}_{3,3,1}}}{p^{12}}} 
 + {\displaystyle \frac {(p^{2}-p)\,{\mathit{sym}_{4,2,1}}}{p^{12}}} \\
 + {\displaystyle \frac {{\mathit{sym}_{4,3,0}}}{p^{10}}}$
\item[15)]
 $\omega ({t}(1,p^{4},p^{4}))
 = {\displaystyle \frac {(p^{2} -2\,p + 1)\,{\mathit{sym}_{3,3,2}}}{p^{14}}} 
 + {\displaystyle \frac {(p^{2}-p)\,{\mathit{sym}_{4,2,2}}}{p^{14}}} 
 + {\displaystyle \frac {(p^{2}-p)\,{\mathit{sym}_{4,3,1}}}{p^{14}}} \\
 + {\displaystyle \frac {{\mathit{sym}_{4,4,0}}}{p^{12}}}$
\item[16)]
 $\omega (t(1,1,p^{5}))
 = {\displaystyle \frac { p^{2} -2\,p+ 1)\,{\mathit{sym}_{2,2,1}}}{p^{7}}} 
 + {\displaystyle \frac {(p^{2} -2\,p+ 1)\,{\mathit{sym}_{3,1,1}}}{p^{7}}} 
 + {\displaystyle \frac {(p^{2}-p)\,{\mathit{sym}_{3,2,0}}}{p^{7}}} \\
 + {\displaystyle \frac {(p^{2}-p)\,{\mathit{sym}_{4,1,0}}}{p^{7}}} 
 + {\displaystyle \frac {{\mathit{sym}_{5,0,0}}}{p^{5}}}$
\item[17)]
 $\omega ({t}(1,p,p^{5}))
 = {\displaystyle \frac {(2\,p^2 - 4\,p + 2)\,{\mathit{sym}_{2,2,2}}}{p^{9}}} 
 + {\displaystyle \frac {(2\,p^{2} - 3\,p+1)\,{\mathit{sym}_{3,2,1}}}{p^{9}}}  
 + {\displaystyle \frac {(p^{2}-p)\,{\mathit{sym}_{3,3,0}}}{p^{9}}} \\ 
 + {\displaystyle \frac {(2\,p^{2} - 2\,p)\,{\mathit{sym}_{4,1,1}}}{p^{9}}} 
 + {\displaystyle \frac {(p^{2}-p)\,{\mathit{sym}_{4,2,0}}}{p^{9}}} 
 + {\displaystyle \frac {{\mathit{sym}_{5,1,0}}}{p^{7}}}$
\item[18)]
 $\omega ({t}(1,p^{2},p^{5}))
 = {\displaystyle \frac {(3\,p^3 - 5\,p^{2} + 3\,p - 1)\,{\mathit{sym}_{3,2,2}}}{p^{12}}} 
 + {\displaystyle \frac {(2\,p^3 - 4\,p^{2} + 2\,p)\,{\mathit{sym}_{3,3,1}}}{p^{12}}}  \\
 + {\displaystyle \frac {(2\,p^3 - 3\,p^{2} + p)\,{\mathit{sym}_{4,2,1}}}{p^{12}}} 
 + {\displaystyle \frac {(p^{3} - p^{2})\,{\mathit{sym}_{4,3,0}}}{p^{12}}} 
 + {\displaystyle \frac {(p^{3} - p^{2})\,{\mathit{sym}_{5,1,1}}}{p^{12}}}  
 + {\displaystyle \frac {{\mathit{sym}_{5,2,0}}}{p^{9}}}$
\item[19)]
 $\omega ({t}(1,p^{3},p^{5}))
 = {\displaystyle \frac {(3\,p^3 - 5\,p^{2} + 3\,p - 1)\,{\mathit{sym}_{3,3,2}}}{p^{14}}} 
 + {\displaystyle \frac {(2\,p^3 - 4\,p^{2} + 2\,p)\,{\mathit{sym}_{4,2,2}}}{p^{14}}} \\  
 + {\displaystyle \frac {(2\,p^3 - 3\,p^{2} + p)\,{\mathit{sym}_{4,3,1}}}{p^{14}}} 
 + {\displaystyle \frac {(p^{3} - p^{2})\,{\mathit{sym}_{4,4,0}}}{p^{14}}} 
 + {\displaystyle \frac {(p^{3} - p^{2})\,{\mathit{sym}_{5,2,1}}}{p^{14}}} 
 + {\displaystyle \frac {{\mathit{sym}_{5,3,0}}}{p^{11}}}$
\item[20)]
 $\omega ({t}(1,p^{4},p^{5}))
 = {\displaystyle \frac {(2\,p^2 - 4\,p + 2)\,{\mathit{sym}_{3,3,3}}}{p^{15}}} 
 + {\displaystyle \frac {(2\,p^2 - 3\,p + 1)\,{\mathit{sym}_{4,3,2}}}{p^{15}}} \\  
 + {\displaystyle \frac {(2\,p^2 - 2\,p    )\,{\mathit{sym}_{4,4,1}}}{p^{15}}} 
 + {\displaystyle \frac {(   p^2 -    p    )\,{\mathit{sym}_{5,2,2}}}{p^{15}}} 
 + {\displaystyle \frac {(   p^2 -    p    )\,{\mathit{sym}_{5,3,1}}}{p^{15}}} 
 + {\displaystyle \frac {                     {\mathit{sym}_{5,4,0}}}{p^{13}}}$
\item[21)]
 $\omega ({t}(1,p^{5},p^{5}))
 = {\displaystyle \frac {(p^2 - 2\,p + 1)\,{\mathit{sym}_{4,3,3}}}{p^{17}}} 
 + {\displaystyle \frac {(p^2 - 2\,p + 1)\,{\mathit{sym}_{4,4,2}}}{p^{17}}} 
 + {\displaystyle \frac {(p^2 -    p    )\,{\mathit{sym}_{5,3,2}}}{p^{17}}} \\ 
 + {\displaystyle \frac {(p^2 -    p    )\,{\mathit{sym}_{5,4,1}}}{p^{17}}} 
 + {\displaystyle \frac {                  {\mathit{sym}_{5,5,0}}}{p^{15}}}$
\item[22)]
 $\omega ({t}(1,1,p^{6}))
 = {\displaystyle \frac {(p^2-2\,p+1)\,{\mathit{sym}_{2,2,2}}}{p^{8}}} 
 + {\displaystyle \frac {(p^2-2\,p+1)\,{\mathit{sym}_{3,2,1}}}{p^{8}}}  
 + {\displaystyle \frac {(p^2-   p  )\,{\mathit{sym}_{3,3,0}}}{p^{8}}} \\ 
 + {\displaystyle \frac {(p^2-2\,p+1)\,{\mathit{sym}_{4,1,1}}}{p^{8}}} 
 + {\displaystyle \frac {(p^2-   p  )\,{\mathit{sym}_{4,2,0}}}{p^{8}}}  
 + {\displaystyle \frac {(p^2-   p  )\,{\mathit{sym}_{5,1,0}}}{p^{8}}} 
 + {\displaystyle \frac {              {\mathit{sym}_{6,0,0}}}{p^{6}}}$
\item[23)]
 $\omega ({t}(1,p,p^{6}))
 = {\displaystyle \frac {(2\,p^2 - 4\,p + 2)\,{\mathit{sym}_{3,2,2}}}{p^{10}}} 
 + {\displaystyle \frac {(2\,p^2 - 3\,p + 1)\,{\mathit{sym}_{3,3,1}}}{p^{10}}} \\ 
 + {\displaystyle \frac {(2\,p^2 - 3\,p + 1)\,{\mathit{sym}_{4,2,1}}}{p^{10}}} 
 + {\displaystyle \frac {(   p^2 -    p    )\,{\mathit{sym}_{4,3,0}}}{p^{10}}} 
 + {\displaystyle \frac {(2\,p^2 - 2\,p    )\,{\mathit{sym}_{5,1,1}}}{p^{10}}} \\
 + {\displaystyle \frac {(   p^2 -    p    )\,{\mathit{sym}_{5,2,0}}}{p^{10}}} 
 + {\displaystyle \frac {                     {\mathit{sym}_{6,1,0}}}{p^{8}}}$
\item[24)]
 $\omega ({t}(1,p^{2},p^{6}))
 = {\displaystyle \frac {(3\,p^3 - 6\,p^2 + 4\,p - 1)\,{\mathit{sym}_{3,3,2}}}{p^{13}}}  
 + {\displaystyle \frac {(3\,p^3 - 5\,p^2 + 3\,p - 1)\,{\mathit{sym}_{4,2,2}}}{p^{13}}} \\
 + {\displaystyle \frac {(2\,p^3 - 4\,p^2 + 2\,p    )\,{\mathit{sym}_{4,3,1}}}{p^{13}}}  
 + {\displaystyle \frac {(   p^3 -    p^2           )\,{\mathit{sym}_{4,4,0}}}{p^{13}}} 
 + {\displaystyle \frac {(2\,p^3 - 3\,p^2 +    p    )\,{\mathit{sym}_{5,2,1}}}{p^{13}}} \\
 + {\displaystyle \frac {(   p^3 -    p^2           )\,{\mathit{sym}_{5,3,0}}}{p^{13}}} 
 + {\displaystyle \frac {(   p^3 -    p^2           )\,{\mathit{sym}_{6,1,1}}}{p^{13}}} 
 + {\displaystyle \frac {                              {\mathit{sym}_{6,2,0}}}{p^{10}}}$
\item[25)]
 $\omega ({t}(1,p^{3},p^{6}))
 = {\displaystyle \frac {(4\,p^3 - 7\,p^2 + 5\,p - 2)\,{\mathit{sym}_{3,3,3}}}{p^{15}}}  
 + {\displaystyle \frac {(3\,p^3 - 6\,p^2 + 4\,p - 1)\,{\mathit{sym}_{4,3,2}}}{p^{15}}} \\
 + {\displaystyle \frac {(2\,p^3 - 4\,p^2 + 2\,p    )\,{\mathit{sym}_{4,4,1}}}{p^{15}}}  
 + {\displaystyle \frac {(2\,p^3 - 4\,p^2 + 2\,p    )\,{\mathit{sym}_{5,2,2}}}{p^{15}}} 
 + {\displaystyle \frac {(2\,p^3 - 3\,p^2 +    p    )\,{\mathit{sym}_{5,3,1}}}{p^{15}}} \\
 + {\displaystyle \frac {(   p^3 -    p^2           )\,{\mathit{sym}_{5,4,0}}}{p^{15}}}  
 + {\displaystyle \frac {(   p^3 -    p^2           )\,{\mathit{sym}_{6,2,1}}}{p^{15}}} 
 + {\displaystyle \frac {                              {\mathit{sym}_{6,3,0}}}{p^{12}}}$
\item[26)]
 $\omega ({t}(1,p^{4},p^{6}))
 = {\displaystyle \frac {(3\,p^3 - 6\,p^2 + 4\,p - 1)\,{\mathit{sym}_{4,3,3}}}{p^{17}}}  
 + {\displaystyle \frac {(3\,p^3 - 5\,p^2 + 3\,p - 1)\,{\mathit{sym}_{4,4,2}}}{p^{17}}} \\
 + {\displaystyle \frac {(2\,p^3 - 4\,p^2 + 2\,p    )\,{\mathit{sym}_{5,3,2}}}{p^{17}}}  
 + {\displaystyle \frac {(2\,p^3 - 3\,p^2 +    p    )\,{\mathit{sym}_{5,4,1}}}{p^{17}}} 
 + {\displaystyle \frac {(   p^3 -    p^2           )\,{\mathit{sym}_{5,5,0}}}{p^{17}}} \\
 + {\displaystyle \frac {(   p^3 -    p^2           )\,{\mathit{sym}_{6,2,2}}}{p^{17}}} 
 + {\displaystyle \frac {(   p^3 -    p^2           )\,{\mathit{sym}_{6,3,1}}}{p^{17}}} 
 + {\displaystyle \frac {                              {\mathit{sym}_{6,4,0}}}{p^{14}}}$
\item[27)]
 $\omega ({t}(1,p^{5},p^{6}))
 = {\displaystyle \frac {(2\,p^2 - 4\,p + 2)\,{\mathit{sym}_{4,4,3}}}{p^{18}}} 
 + {\displaystyle \frac {(2\,p^2 - 3\,p + 1)\,{\mathit{sym}_{5,3,3}}}{p^{18}}} \\ 
 + {\displaystyle \frac {(2\,p^2 - 3\,p + 1)\,{\mathit{sym}_{5,4,2}}}{p^{18}}} 
 + {\displaystyle \frac {(2\,p^2 - 2\,p    )\,{\mathit{sym}_{5,5,1}}}{p^{18}}} 
 + {\displaystyle \frac {(   p^2 -    p    )\,{\mathit{sym}_{6,3,2}}}{p^{18}}} \\
 + {\displaystyle \frac {(   p^2 -    p    )\,{\mathit{sym}_{6,4,1}}}{p^{18}}} 
 + {\displaystyle \frac {                     {\mathit{sym}_{6,5,0}}}{p^{16}}}$
\item[28)]
 $\omega ({t}(1,p^{6},p^{6}))
 = {\displaystyle \frac {(p^2 - 2\,p + 1)\,{\mathit{sym}_{4,4,4}}}{p^{20}}} 
 + {\displaystyle \frac {(p^2 - 2\,p + 1)\,{\mathit{sym}_{5,4,3}}}{p^{20}}} \\
 + {\displaystyle \frac {(p^2 - 2\,p + 1)\,{\mathit{sym}_{5,5,2}}}{p^{20}}} 
 + {\displaystyle \frac {(p^2 -    p    )\,{\mathit{sym}_{6,3,3}}}{p^{20}}} 
 + {\displaystyle \frac {(p^2 -    p    )\,{\mathit{sym}_{6,4,2}}}{p^{20}}} \\
 + {\displaystyle \frac {(p^2 -    p    )\,{\mathit{sym}_{6,5,1}}}{p^{20}}} 
 + {\displaystyle \frac {                  {\mathit{sym}_{6,6,0}}}{p^{18}}}$ \, . \qed
\end{itemize}

Note that Rhodes, J. A. and Shemanske, T. R. developed an alternative method of computing
$\omega(t(p^{\delta_1},\dots ,p^{\delta_n}))$ in \cite{RhSh}, based on counting of certain left
cosets in a given double coset (Theorem 4.3 of \cite{RhSh}). 

\section{A special case}
For some particular values of the Satake parameters $x_0$, $x_1$, $x_2$, $x_3$, the polynomial $P_3$
can be considerably simplified.  For example,
let us substitute $x_0=1$, $x_1=p$, $x_2=p^2$ and $x_3=p^3$ as in Exercise 3.3.40, p.168 of \cite{An87}:
$\mathit{P}_\nu (v) :=\mathit{P}(1, p, p^2, p^3, v)$, where $\nu$ denotes the degree homomorphism
$\nu(x_0)=1$, $\nu(x_1)=p$, $\nu(x_2)=p^2$, $\nu(x_3)=p^3$. Then the polynomial $P$ takes the form

\begin{align*}&
\mathit{P}_\nu (v) = 1 - \left({\displaystyle \frac {p^{7} + p^{8} + p^{9
}}{p}} + (p^{2} + p + 1)\,p^{4} + {\displaystyle \frac {p^{3} + 
p^{4} + p^{5}}{p}} \right)\,v^{2} 
\\ &
\mbox{} + \left((p + 1)\,p^{10} + {\displaystyle \frac {(p + 1)\,(p^{9
} + p^{10} + p^{11})}{p^{2}}}  + {\displaystyle \frac {(p + 1)\,(
p^{7} + p^{8} + p^{9})}{p^{2}}}  + (p + 1)\,p^{4}\right)\,v^{3} \\ &
\mbox{} - \left({\displaystyle \frac {p^{13} + p^{14} + p^{15}}{p^{
2}}} + (p^{2} + p + 1)\,p^{9}+ {\displaystyle \frac {p^{9} + p
^{10} + p^{11}}{p^{2}}} \right)\,v^{4} + p^{15}\,v^{6} \\ &
=1 - (p^{8}\,+ p^{7}\,+2\,p^{6}\,+ p^{5}\,+2\,p^{4}\,+ p^{3}+p^{2})v^{2}  \\ &
+(p^{11}+ 2\,p^{10} + 2\,p^{9}+ 3\,p^{8}   + 3\,p^{7} + 2\,p^{6} + 
2\,p^{5} + p^{4})\,v^{3} \\ &
-(p^{13}\,+ p^{12}\,+2\,p^{11}\,+ p^{10}\,+2\,p^{9}\,+ p^{8}+p^{7})v^{4}
+ p^{15}\,v
^{6}. 
\end{align*}
This gives the following factorization:
\begin{align*}
\mathit{P}_\nu (v)
=(1-p^{}\,v)\,(1-p^{2}\,v )\,(1-p^{3}\,v )\,(1-p^{4}\,v)\,(1+ p^{}\,v + p^{2}\,v + p^{3}\,v + p^{4}\,v +p^{5}\,v^{2} ).
\end{align*}
%

\bibliographystyle{amsplain}

\end{document}